\documentstyle[11pt,leqno]{article}
\newtheorem{theorem}{{\sc Theorem}}
\newcommand{\bt}{\begin{theorem}}
\newcommand{\et}{\end{theorem}}
\newcommand{\newsection}[1]{\setcounter{equation}{0} \setcounter{theorem}{0}
\section{#1}}

\newcommand{\NI}{\noindent}
\newcommand{\bea}{\begin{eqnarray}}
\newcommand{\eea}{\end{eqnarray}}

\def \spec#1 {\mathop{#1}}

\def \b #1 {\bf #1}

\newcommand {\CC}{\centerline}

\newcommand{\clx}{{\cal X}}

\newcommand{\ity}{\infty}
\newcommand{\raro}{\rightarrow}

\newcommand{\vsp}{\vskip 1em}

\newcommand{\be}{\begin{equation}}
\newcommand{\ee}{\end{equation}}
\newcommand{\ben}{\begin{eqnarray*}}
\newcommand{\een}{\end{eqnarray*}}

\begin{document}
\CC {\bf  }
\CC{\bf CHARACTERIZATION OF PROBABILITY MEASURES BASED ON }
\CC{\bf Q-INDEPENDENT GENERALIZED RANDOM FIELDS}
\vsp
\CC { B.L.S. Prakasa  Rao{\footnote {E-mail address: blsprao@gmail.com}}}
\CC{  CR Rao Advanced Institute of  Mathematics, Statistics and }
\CC{ Computer  Science, Hyderabad 500046, India}
\vsp
\NI{\bf Abstract:} Prakasa Rao ({\it Studia Sci. Math. Hungar.}, {\bf 11} (1976) 277-282) studied a characterization of probability distributions  for linear functions of independent generalized random fields. These results are extended to $Q$-independent generalized random fields. It is known that independence of random variables implies $Q$-independence of them but the converse is not true.
\vsp
\NI{\bf Key words :} Generalized random field; Characteristic functional; $Q$-independence; Gaussian characteristic functional. 
\vsp
\NI{\bf MSC2020:} Primary 60G60. 
\vsp
\newsection{Introduction} Rao (1971) obtained some characterizations of probability distributions on the real line through linear functions of independent real-valued random variables. Some of these results were extended to linear functions of independent generalized random fields in Prakasa Rao (1976). Kagan and Szekely (2016) introduced the concept of $Q$-independence for real-valued random variables. A characterization of probability distributions for $Q$-independent random elements was presented in a collection of articles in Prakasa Rao (2016, 2017, 2018a,b,c) and for $Q$-independent random variables taking values in a locally compact Abelian group  by Feldman (2017). We now extend the results in Prakasa Rao (1976) to $Q$-independent generalized random fields. It is known that independence of random variables implies their $Q$-independence but the converse is not true.
\vsp 
\newsection{Preliminaries} Let $\clx$ be the space of all real-valued functions $\phi(x)= \phi(x_1,\dots,x_n)$ of $n$ variables which are infinitely differentiable and have bounded supports. A sequence of functions $\{\phi_m, m \geq 1\}$ in $\clx$ is said to converge to zero if there exists a constant $a$ such that $\phi_m$ vanishes for $||x||\geq a,$ and if , for every integer $q\geq 1,$ the sequence $\{\phi_m^{(q)}, m \geq 1\}$ converges uniformly to zero where $||x||= (x_1^2+\dots+x_n^2)^{1/2}$ and $\phi^{(q)}$ denotes the $q$-th derivative of $\phi.$Any continuous linear functional on $\clx$ is called a {\it generalized function}.
\vsp
A functional $\Phi$ defined on $\clx$ is said to be a {\it random functional} if for every $\phi\in \clx$ there is associated a real-valued random variable $\Phi(\phi).$ In other words, for every set of $m$ elements $\phi_i, 1 \leq i \leq m$ in $\clx$ , one can specify the probability that 
$$a_i\leq \Phi(\phi_i)\leq b_i, 1 \leq i \leq m$$
for $-\ity <a_i<b_i<\ity, 1 \leq i \leq m$ and these probability distributions are consistent. The random fuctional $\Phi$ is said to be {\it linear} if for any two elements $\phi, \psi \in \clx,$ and for any two real numbers $\alpha, \beta,$ 
$$\Phi(\alpha \phi+\beta \psi)= \alpha \Phi(\phi) +\beta \Phi(\psi)$$
almost surely. A random functional $\Phi$ is said to be {\it continuous} if the convergence of the functions $\phi_{kj}$ to $\phi_j, 1\leq j \leq m$ as $k \raro \ity$ in $\clx$ implies that for every bounded continuous function $f(x_1,\dots, x_m),$
$$\lim_{k \raro \ity}\int_{R^m} f(x_1,\dots,x_m)P_k(dx)=\int_{R^m} f(x_1,\dots,x_m)P(dx)$$
where $P$ is the probability measure corresponding to the random vector $(\Phi(\phi_1),\dots,\Phi(\phi_m))$ and $P_k$ is the probability measure corresponding to the random vector $(\Phi(\phi_{k1}),\dots,\Phi(\phi_{km})).$
\vsp
Any continuous linear random functional on $\clx$ is called a {\it generalized random function}. If the space $\clx$ consists of functions of one variable, then the corresponding random functional is called a {\it generalized random process}. If the space $\clx$ consists of functions of several variables, then the corresponding random functional is called a {\it generalized random field}.
\vsp
Let $\Phi$ and $\Psi$ be two generalized random fields on $\clx.$ The generalized random fields $\Phi$ and $\Psi$ are said to be {\it independent} if the set of random variables $\{\Phi(\phi), \phi\in \clx\}$  is independent of  the set of random variables $\{\Psi(\phi), \phi\in \clx\}.$ This notion can be extended to any finite number of generalized random fields in an obvious manner.
\vsp
Let $\Phi$ be a generalized random field. The functional
$$L(\phi)= E[e^{i\Phi(\phi)}], \phi \in \clx$$
is called the {\it characteristic functional} of the generalized random field $\Phi.$ It can be shown that $L(0)=1, L(-\phi)= L(\bar \phi), L(\phi)$ is continuous in $\phi$ and positive definite. Conversely, if $L(.)$ is a positive-definite continuous functional on $\clx$ such that $L(0)=1,$ it can be shown that there exists a generalized random field $\Phi$ on $\clx$ whose characteristic functional is $L(.)$. Furthermore the correspondence between the characteristic functionals $L(.)$ and the generalized random fields $\Phi$ on $\clx$ is one to one.
\vsp
Let $\Phi_1,\dots,\Phi_k$ be generalized random fields on the space $\clx.$ The joint characteristic functional of the $k$-dimensional generalized random field $(\Phi_1,\dots,\Phi_k)$ is defined by
$$L_{\Phi_1,\dots,\Phi_k}(\phi_1,\dots,\phi_k)= E[\exp\{i \Phi_1(\phi_1)+\dots+i\Phi_k(\phi_k)\}], \phi_j, 1 \leq j \leq k \in \clx.$$ 
\vsp
If the generalized random fields are independent, then it can be shown that
$$L_{\Phi_1,\dots,\Phi_k}(\phi_1,\dots,\phi_k)=L_1(\phi_1)\dots L_k(\phi_k)$$
where $L_j(.)$ is the characteristic functional of $\Phi_j$ for $1\leq j \leq k.$ 
\vsp
A generalized random field $\Phi$ on $\clx$ is said to be {\it Gaussian} if its characteristic functional is of the form
$$L(\phi)= \exp(i\; m(\phi)- \frac{1}{2}B(\phi,\phi)), \phi \in \clx$$
where $m(.)$ is a generalized function and $B(\phi,\psi)= E[\Phi(\phi)\Phi(\psi)], \phi,\psi \in \clx.$ It is said to be {\it degenerate} if 
its characteristic functional is of the form
$$L(\phi)= \exp (i\; m(\phi)), \phi \in \clx$$
where $m(.)$ is a generalized function.
\vsp
We refer the reader to Gelfand and Vilenkin (1964) for more details on generalized random fields and generalized random processes.
\vsp
Denote by $\Delta_h$ the finite difference operator
$$\Delta_h f(\phi)= f(\phi+ h)-f(h). $$
a function $f(\phi)$ defined on $\clx$ is called a {\it polynomial} if 
$$\Delta_h^{n+1}f(\phi)=0$$
for some integer $n \geq 1$ and for all $\phi, h \in \clx.$ The minimal integer $n$ for which this equality holds is called the {\it degree of the polynomial} $f(.)$ defined on $\clx.$
\vsp
Let $\Phi_i, 1\leq i \leq k$ be generalized random fields on $\clx.$ Let $\beta_i, 1 \leq i \leq k$ be nonzero real numbers. We define the process $\beta_1\Phi_1+\dots+\beta_k\Phi_k$ to be the process for which to every $\phi \in \clx$ corresponds the random variable $\Phi_1(\beta_1\phi)+\dots+\Phi_k(\beta_k\phi).$
\vsp
Let $\Phi_1,\dots,\Phi_k$ be generalized random fields. We say that they are $Q$-{\it independent} if their joint characteristic functional 
can be represented in the form
\be
L_{\Phi_1,\dots,\Phi_k}(\phi_1,\dots,\phi_k)= \Pi_{i=1}^k L_{\Phi_i}(\phi_i) \exp(q(\phi_1,\dots,\phi_k)), \phi_1, 1\leq i \leq k \in \clx
\ee
where $q(\phi_1,\dots,\phi_k)$ is a continuous  polynomial on the space $\clx^k$ and $q(0,\dots,0)=0.$
\vsp
Suppose that $\Phi_i, i=1,2,3$ are independent Gaussian random fields. Then it is obvious that $\eta_1=\Phi_1+ \Phi_2$ and $\eta_2= \Phi_1+\Phi_3$ are not independent random fields. However they are $Q$-independent. This can be seen by computing the joint characteristic functional of the bivariate generalized random field $(\eta_1,\eta_2)$ and the characteristic functionals of the generalized random fields $\eta_1$ and $\eta_2.$ 
\vsp
The following result is a consequence  of the Marcinkeiwicz theorem (cf. Marcinkiewicz (1938)) for real-valued random variables.
\vsp
\NI{\bf Theorem 2.1:} Let $f(y)$ be the characteristic functional of a  generalized random field $\Phi$ on $\clx.$ If 
$$ f(y)=\exp[P(y)] , y \in \clx,$$
where $P(y)$ is a continuous polynomial in $y \in \clx,$  then $P(y)$ is a polynomial of degree less than or equal to 2 and $f(y)$ is the characteristic functional of a Gaussian random field which could be degenerate.
\vsp
\NI{\bf Proof :} By the definition of the characteristic functional of the  generalized random field $\Phi$ on $\clx,$ it follows that
$$E[\exp(i \Phi(y))]= \exp[P(y)], y \in \clx.$$
Hence
$$E[ \exp (i t \Phi(y))]= E[\exp(i\Phi(ty))] = \exp[P(ty)]$$
for any real $t \in R.$ Since the function on the left side of the above equation is the characteristic function of the real valued random variable $\Phi(y)$ it follows that the function on the right side of the equation has to be a polynomial of degree less than or equal to 2 in $ty$ by the classical Marcinkeiwicz theorem. Choosing $t=1,$ it follows that $P(y)$ is a polynomial of degree less than or equal to 2 in $y$ which in turn implies that the generalized random field is either degenerate or is Gaussian. 
\vsp
We now prove a theorem dealing with functional equations on the space $\clx $ which is of independent interest.
Proof of the theorem is similar to that when the space $\clx$ is the set of real numbers (cf. Kagan et al. (1973)). Our presentation is similar to that in Feldman (2017) when the space $\clx$ is a locally compact Abelian group. We present the detailed proof for completeness.
\vsp
\NI{\bf Theorem 2.2 :} Let $\clx $ be the space of infinitely differentiable functions. Consider the functional equation
\be
\sum_{j=1}^n \psi_j(u+b_jv)= P(u)+ Q(v)+R(u,v), u,v \in \clx 
\ee
where $b_1,\dots,b_n$ are nonzero real numbers with $b_i \neq b_j,1\leq j \leq n$ and  $\psi_j(u),1 \leq j \leq n, P(u), Q(v)$ are functions on $\clx $ and $R(u,v)$ is a polynomial on $\clx \times \clx.$ Then $P(u)$ is a polynomial on $\clx.$
\vsp
\NI{\bf Proof:} We use the finite difference method for proving the theorem. Let $h_1$ be an arbitrary element of $\clx.$ Define $k_1= -b_n^{-1}h_1.$ Then $h_1+b_n k_1=0.$ Substitute $u+h_1$ for  $u$ and $v+k_1$ for $v$ in the equation (2.2). Subtracting the equation (2.2) from the resulting equation, it follows that
\be
\sum_{j=1}^{n-1}\Delta_{\ell_{1j}} \psi_j(u+b_j v)= \Delta_{h_1}P(u)+ \Delta_{k_1}Q(v)+ \Delta_{(h_1,k_1)} R(u,v), u,v \in \clx 
\ee
where $\ell_{1j}= h_1+b_j k_1= (b_j-b_n)k_1, j=1,\dots, n-1.$ Let $h_2$ be an arbitrary element of $\clx.$ Let $k_2= -b_{n-1}^{-1}h_2.$ Then $h_2+b_{n-1}k_2=0.$ Substitute $u+h_2$ for $u$ and $v+k_2$ for $v$ in the equation (2.3). Subtracting equation (2.3) from the resulting equation, it follows that
\bea
\sum_{j=1}^{n-2}\Delta_{\ell_{2j}}\Delta_{\ell_{1j}}\psi_j(u+b_j v) &= & \Delta_{h_2}\Delta_{h_1}P(u)+ \Delta_{k_2}\Delta_{k_1}Q(v)\\\nonumber
&&\;\;\;\;+\Delta_{(h_2,k_2)}\Delta_{(h_1,k_1)}R(u,v), u, v \in \clx, \\\nonumber
\eea
where $\ell_{2j}= h_2+b_jk_2=(b_j-b_{n-1})k_2, j=1,\dots,n-2.$ Following similar arguments. we get the equation
\bea
\;\;\;\;\\\nonumber
\Delta_{\ell_{n-1,1}}\Delta_{\ell_{n-2,1}}\dots \Delta_{\ell_{1,1}}\psi_1(u+b_1v) &=& \Delta_{h_{n-1}}\Delta_{h_{n-2}}\dots\Delta_{h_1}P(u)\\\nonumber
&&\;\;\;\;+\Delta_{k_{n-1}}\Delta_{k_{n-2}}\dots\Delta_{k_1}Q(v)\\\nonumber
&&\;\;\;\;+ \Delta_{(h_{n-1}, k_{n-1})}\Delta_{(h_{n-2}, k_{n-2})}\dots\Delta_{(h_1, k_1)}R(u,v),  \\\nonumber
\eea
for $u,v \in \clx$, where $h_m$ are arbitrary elements in $\clx, k_m= -b_{n-m+1}^{-1}h_m, m=1,2,\dots,n-1, \ell_{mj}= h_m+b_jk_m= (b_j-b_{n-m+1}k_m, j=1,2,\dots n-m.$ Let $h_n$ be an arbitrary element of $\clx.$ Let $k_n=-b_1^{-1}h_n.$ Then  $h_n+b_1k_n=0.$ Substitute $u+h_n$ for $u$ and $v+k_n$ for $v$ in the equation (2.5). Subtracting the equation (2.5) from the resulting equation, we get that
\bea
\Delta_{h_n}\Delta_{h_{n-1}}\dots\Delta_{h_1}P(u)+ \Delta_{k_n}\Delta_{k_{n-1}}\dots\Delta_{k_1}Q(v)\\\nonumber
\;\;\;\;+\Delta_{(h_n,k_n)}\Delta_{(h_{n-1},k_{n-1})}\dots\Delta_{(h_1,k_1)}R(u,v)=0, u,v \in \clx .\\\nonumber
\eea
Let $h_{n+1}$ be an arbitrary element of $\clx.$ Substitute $h_{n+1}$ for $u$ in the equation (2.6). Subtracting the equation (2.6) from the resulting equation, we obtain that
\bea
\Delta_{h_{n+1}}\Delta_{h_n}\Delta_{h_{n-1}}\dots\Delta_{h_1}P(u)\\\nonumber
\;\;\;\;+\Delta_{(h_n,k_n)}\Delta_{(h_{n-1},k_{n-1})}\dots\Delta_{(h_1,k_1)}R(u,v)=0, u,v \in \clx .\\\nonumber
\eea
Observe that, if $h$ and $k$ are arbitrary elements of the space $\clx,$ it follows that
\be
\Delta_{(h,k)}^{\ell+1}R(u,v)=0, u,v \in \clx 
\ee
for some integer $\ell \geq 0$ since $R(u,v)$ is a polynomial in $(u,v)$ by hypothesis. Since $h_m, m=1,\dots, n+1$ are arbitrary elements of the space $\clx,$ let us choose $h_1=\dots = h_{n+1}=h \in \clx$ in the equation (2.7) and apply the operator $\Delta_{(h,k)}^{\ell+1}$ to both sides of the resulting equation. Applying the equation (2.8) now leads to the equation
\be
\Delta_h^{\ell+n+2}P(u)=0, u,h \in \clx.
\ee
Hence the function $P(u)$ is a polynomial of degree at most $\ell+n+1.$
\vsp
\NI{\bf Remarks:} Let $\ell$ be the degree of the polynomial $R(u,v)$ in Theorem 2.2. Following the methods in Kagan et al. (1973), it can be shown that the degree of the polynomial $P(u)$ in Theorem 2.1 does not exceed $\max(n,\ell)$ where $n$ is the number of functions in the left side of the functional equation (2.2).
\vsp
Two generalized random fields $\Phi$ and $\Psi$ are said to be ``determined up to a Gaussian generalized random field" if there exist a  generalized random field $\Lambda$ such that $\Phi= \Psi+\Lambda$ almost surely. They are said to be determined up to ``translation" if there exists a generalized function $m$ such that $\Phi=\Psi+m$ almost surely.
\vsp
\newsection{Main Results} 
We now prove a theorem characterizing generalized random fields up to Gaussian factors. 
\vsp
\NI{\bf Theorem 3.1:} Let $\Phi_i, 0 \leq i \leq 3$ be four $Q$-independent generalized random fields on $\clx$ and let
\bea
\Psi_1 &=& \Phi_0+\Phi_1+\Phi_2+\Phi_3\\\nonumber
\Psi_2 &=& \beta_0\Phi_0+\beta_1\Phi_1+\beta_2\Phi_2+\beta_3\Phi_3\\\nonumber
\eea
where $\beta_i, 0\leq i \leq 3$ are non-zero real numbers such that $\beta_i\neq \beta_j, 0 \leq i \neq j \leq 3.$ Further suppose that the joint characteristic functional $H(\phi,\psi)$ of $(\Psi_1,\Psi_2)$ does not vanish. If $L_i(\phi)$ and $M_i(\phi)$ are two alternate possible characteristic functionals of the generalized random field $\Phi_i, 0\leq i \leq 3,$ then
\be
L_j(\phi)= M_j(\phi)\exp(i\;m_j(\phi)-\frac{1}{2}B_j(\phi,\phi)), 0\leq j \leq 3
\ee
for some generalized functions $m_j(\phi), 0\leq j \leq 3$ and for some continuous bilinear Hermitian functionals $B_j(\phi,\psi), 0\leq j \leq 3.$
\vsp
\NI{\bf Proof:} Let $\Gamma_i, 0\leq i \leq 3$ be $Q$-independent generalized random fields on $\clx$ such that the two-dimensional generalized random field $(\Sigma_1,\Sigma_2)$ where
\bea
\Sigma_1 &= & \Gamma_0+\Gamma_1+\Gamma_2+\Gamma_3\\\nonumber
\Sigma_2 &= & \beta_0\Gamma_0+\beta_1\Gamma_1+\beta_2\Gamma_2+\beta_3\Gamma_3\\\nonumber
\eea
has the same joint characteristic functional $H(\phi,\psi)$ as that of $(\Psi_1,\Psi_2).$ Let $L_i(.)$ and $M_i(.), 0\leq i \leq 3$ be the characteristic functionals of $\Phi_i$ and $\Gamma_i, 0\leq i \leq 3$ respectively. From the $Q$-independence of the generalized random fields $\Phi_i, 0 \leq i \leq 3,$ it follows that
$$H(\phi,\psi)= \Pi_{i=0}^3M_i(\phi+\beta_i\psi) \exp(P_1(\phi,\psi)), \phi,\psi\in \clx$$
for some polynomial $P_1(\phi,\psi).$ From the $Q$-independence of the generalized random fields $\Gamma_i, 0 \leq i \leq 3,$ it follows that
$$H(\phi,\psi)= \Pi_{i=0}^3L_i(\phi+\beta_i\psi) \exp(P_2(\phi,\psi)), \phi,\psi\in \clx$$
for some polynomial $P_2(\phi,\psi).$ Hence
\bea
H(\phi,\psi) &=& \Pi_{i=0}^3M_i(\phi+\beta_i\psi) \exp(P_1(\phi,\psi))\\\nonumber
&=& \Pi_{i=0}^3L_i(\phi+\beta_i\psi) \exp(P_2(\phi,\psi)), \phi,\psi \in \clx.\\\nonumber
\eea
Since $H(\phi,\psi) \neq 0$ for all $\phi,\psi \in \clx$ by hypothesis, the equation given above implies that $L_i(\phi+\beta_i\psi) \neq 0, 0\leq i \leq 3$ and $M_i(\phi+\beta_i\psi) \neq 0, 0\leq i \leq 3$ for all $\phi,\psi \in \clx.$ Let
$$J_i(\phi)=\log \frac{L_i(\phi)}{M_i(\phi)}, 0\leq i \leq 3$$
where the logarithm is taken to be the continuous branch with $J_i(0)=0.$ The equation (3.4) implies that
\be
\sum_{i=0}^3J_i(\phi+\beta_i\psi)= P_1(\psi,\phi)-P_2(\psi,\phi), \phi, \psi \in \clx  
\ee   
where $P_1(.,.)$ and $P_2(.,.)$ are polynomials. Since $\beta_i \neq \beta_j, 0\leq i\neq j \leq 3$ and $\beta_j \neq 0,$ applying arguments similar to those in the proof of Lemma 1.5.1 in Kagan et al. (1973), it follows that the functions $J_i(\phi), i=0,\dots,3$ are polynomials in $\phi$ on $\clx.$  Hence there exists polynomials $f_j(\phi)$ such that
\be
L_j(\phi)= M_j(\phi)\exp[f_j(\phi)], \phi \in \clx, 0\leq j \leq 3.
\ee
Note that the functional $L_j(,)$ on the left side of the equation (3.6)  is a characteristic functional and it is non-vanishing by the equation  (3.4). Hence the function on the right side of the equation is also a non-vanishing characteristic functional which in turn implies that the functional $\exp[f_j(\phi)], \phi \in \clx$ is a characteristic functional by the one-to-one correspondence between the probability measures and the characteristic functionals on the space $\clx.$ An application of the  Marcinkeiwicz lemma (cf. Theorem 2.1) implies that the degree of the polynomial $f_j(.)$ can not exceed two. It can be shown that 
\be
L_j(\phi)= M_j(\phi)\exp(i\; m_j(\phi)-\frac{1}{2}B_j(\phi,\phi)), 0\leq j \leq 3
\ee
for some generalized functions $m_j(\phi), 0\leq j \leq 3$ and for some continuous bilinear Hermitian functional $B_j(\phi,\psi), 0\leq j \leq 3$ by arguments similar to those in Prakasa Rao (1976), p.281. 
\vsp
The following theorem can be proved by arguments similar to those given above. We omit the details.
\vsp
\NI{\bf Theorem 3.2:} Let $\Phi_i, 0 \leq i \leq 2$ be four $Q$-independent generalized random fields on $\clx$ and let
\bea
\Psi_1 &=& \Phi_0+\Phi_1+\Phi_2\\\nonumber
\Psi_2 &=& \beta_0\Phi_0+\beta_1\Phi_1+\beta_2\Phi_2\\\nonumber
\eea
where $\beta_i, 0\leq i \leq 2$ are non-zero real numbers such that $\beta_i\neq \beta_j, 0 \leq i \neq j \leq 2.$ Further suppose that the joint characteristic functional $H(\phi,\psi)$ of $(\Psi_1,\Psi_2)$ does not vanish. If $L_i(\phi)$ and $M_i(\phi)$ are two alternate possible characteristic functionals of the generalized random field $\Phi_i, 0\leq i \leq 2,$ then
\be
L_j(\phi)= M_j(\phi)\exp(i\;m_j(\phi)), 0\leq j \leq 2
\ee
for some generalized functions $m_j(\phi), 0\leq j \leq 2.$ 
\vsp
\NI{\bf Acknowledgment :} Work on this paper was supported by the scheme ``INSA Senior Scientist" at the  CR Rao Advanced Institute of Mathematics, Statistics and Computer Science, Hyderabad 500046, India.
\vsp
\NI{\bf References}
\begin{description}
\item Feldman, G. (2017) Characterization theorems for $Q$-independent random variables with values in a locally compact Abelian group, {\it Aequat. Math.}, {\bf 91}, 949-967. 
\item Gelfand, I.M. and Vilenkin, N. Ya. (1964) {\it Generalized Functions}, Vol. 4, Academic Press, New York.
\item Kagan, A.M., Linnik, Yu.V., and Rao, C.R. (1973) {\it Characterization Problems in Mathematical Statistics}, Wiley, New York. 
\item Kagan, A.M. and Szekely, G.J. (2016) An analytic generalization of independence and identical distributiveness, {\it Statist. Probab. Lett.}, {\bf 110}, 244-248.
\item Marcinkiewicz, J. (1938) Sur une propriete de la loi Gauss, {\it Math. Zeit.}, {\bf 44}, 612-618.
\item Prakasa Rao, B.L.S. (1976) On a property of generalized random fields, {\it Studia Sci. Math. Hungar.}, {\bf 11}, 277-282.
\item Prakasa Rao, B.L.S. (2016) Characterization of probability distributions through linear forms of $Q$-conditionally independent random variables, {\it Sankhya Series A}, {\bf 78-A}, 221-230. 
\item Prakasa Rao, B.L.S. (2017) Characterization of probability measures on Hilbert spaces via $Q$-independence, {\it J. Indian Stat. Assoc.}, {\bf 55}, 95-106.
\item Prakasa Rao, B.L.S. (2018a) Characterization of probability measures on locally compact Abelian groups via $Q$-independence, {\it Acta Math. Szeged}, {\bf 84}, 705-711.
\item Prakasa Rao, B.L.S. (2018b) Characterization of probability distributions through $Q$-independence, {\it Theory Probab. Appl.}, {\bf 62}, 335-338.
\item Prakasa Rao, B.L.S. (2018c) On the Skitovitch-Darmois-Ramachandran-Ibragimov theorem for linear forms of $Q$-independent  random sequences, {\it Studia Sci. Math. Hungar.}, {\bf 55}, 353-363.
\item Rao, C.R. (1971) Characterization of probability laws through linear functions, {\it Sankhya, Series A}, {\bf 33}, 255-259.
\end{description}
\end{document}